\documentclass[11pt]{amsart} 
\usepackage{amsmath}
\usepackage{amssymb}
\usepackage{amscd}
\usepackage{latexsym}
\usepackage{graphics, epsfig}
\usepackage{curves}

\theoremstyle{plain}
\newtheorem{thm}{Theorem}[section]
\newtheorem{lem}[thm]{Lemma}

\theoremstyle{definition}
\newtheorem{defn}[thm]{Definition}

\def \R {\mathbf{R}}
\def \Z {\mathbf{Z}}
\def \Sig{\Sigma}

\def \CM {\mathcal M}

\def \b {\beta}

\def \o {\omega}

\def \- {\setminus}

\def \ssw {\text{SW}}

\def \DD {\Delta}

\def\fs{\mathfrak{s}}

\title{Construction of new symplectic cohomology $S^{2}\times S^{2}$} 
              
\begin{document}

\author{Anar Akhmedov}
\address{Department of Mathematics\\
University of California at Irvine\\
Irvine, CA 92697-3875}

\email{ahmadov@math.gatech.edu}

\date{November 5, 2006}


\begin{abstract} In this article, we present new symplectic 4-manifolds with same integral cohomology as $S^{2}\times S^{2}$. The generalization of this construction is given as well, an infinite family of symplectic 4-manifolds cohomology equivalent to $\#_{(2g-1)}{(S^{2}\times S^{2})}$ for any $g\geq 2$. We also compute the Seiberg-Witten invariants of these manifolds.   
\end{abstract}

\maketitle

\setcounter{section}{0}

\section{Introduction}

The aim of this article is to construct new examples of symplectic $4$-manifolds with the same integral cohomology as $S^{2}\times S^{2}$. The similar problems have been studied in the algebro-geometric category, i.e. the existence of algebraic surfaces homology equivalent but not isomorphic to $\mathbf{P^{2}}$ (or $\mathbf{P^{1}\times P^{1}}$) as an algebraic variety. Mumford \cite {Mu} and Pardini \cite {P} gave the constructions of such fake $\mathbf{P^{2}}$ and $\mathbf{P^{1}\times P^{1}}$.           

 We study this problem in the symplectic category and prove the following results.

\begin {thm} Let $K$ be a genus one fibered knot in $S^{3}$. Then there exist a symplectic $4$-manifold $X_{K}$ cohomology equivalent to $S^{2}\times S^{2}.$ 

\end{thm}

\begin {thm} Let $K'$ be any genus $g$ fibered knot in $S^{3}$. Then there exist an infinite family of symplectic $4$-manifolds $V_{K'}$  that is cohomology equivalent to $\#_{(2g-1)}{(S^{2}\times S^{2})}$ for $g \geq 2$.

\end{thm}

This paper is organized as follows: Section $2$ contains the basic definitions and formulas that will be important throughout this paper. Section $3$ gives quick introduction to Seiberg-Witten invariants. The remaning two sections are devoted to the construction of family of symplectic 4-manifolds cohomology equivalent to $\#_{(2g-1)}{(S^{2}\times S^{2})}$ and the fundamental group computation for our examples.   

{\bf Acknowledgments:} I am very grateful to Ron Stern for his guidance and support. I am also indebted to B. Doug Park and Andr\'as Stipsicz for the remarks and corrections on earlier drafts of this paper.

\section{Preliminaries}

\subsection{Knots and $3$-Manifolds}

In this section, we give short introduction to the fibered knots, Dehn surgery and state a few facts that we will be used in our construction. We refer the reader to Section 10.H \cite {R} for a more complete treatment.

\begin{defn}\emph{(Torus Knots)} The knot which wraps around the solid torus in the longitudinal direction $p$ times and in the meridinal direction $q$ times called $(p,q)$ torus knot and denoted as $T_{p, q}$. \end{defn}

\begin{lem}\cite{S} a) The group of the torus knot $T_{p,q}$ can be represented as follows:
\begin{eqnarray*} \pi_{1}(S^{3}/T_{p,q}) &=& \langle u, v \mid u^{p} = v^{q} \rangle \end{eqnarray*} 
b) The elements $m = u^{s}v^{r}$, $l = u^{p}m^{-pq}$, where $pr + qs = 1$, describe meridian and longitude of the $T_{p,q}$ for a suitable chosen basepoint.

\end{lem}

\noindent All torus knots belongs to the larger category of fibered knots. 
    
\begin{defn} A knot $K$ in $S^{3}$ is called \emph{fibered} if there is fibration $ f : S^{3}\- K \longrightarrow S^{1}$ behaving ``nicely'' near $K$. This means that $K$ has a neighbourhood framed as $S^{1}\times D^{2}$, with $K \cong S^{1}\times 0$   and restriction of the map $f$ to $S^{1}\times (D^{2}- 0)$ is the map to $S^{1}$ given by $(t,x) \rightarrow x/|x|$ 
\end{defn}

It follows from the definition that a preimage for each point $p\in S^{1}$ is Seifert surface for the given knot. The genus of this Seifert surface will be called the genus of the given fibered knot.  

The fibered knots form a large class among the all classical knots. Below we state two theorems that can be used to detect if the given knot is fibered or not.  

\begin{thm} \cite {S} The knot $K \subset S^{3}$ is a fibered knot of genus $g$ if and only if the commutator subgroup of its knot group $\pi_{1}(S^{3}\backslash K)$ is finitely-generated and free group of rank $2g$.  
\end{thm}

\begin{thm} \cite {BZ} The Alexander polynomial $\DD_{K}(t)$ of a fibered knot in $S^{3}$ is monic i.e. the first and last non-zero coefficients of $\DD_{K}(t)$ are $\pm 1$. 
\end{thm}

\noindent If a genus one knot is fibered, then it can be shown by Theorem 1.5 and also by explicit constructon \cite {BZ} that it is either trefoil or figure eight knot. Also, one can construct infinitely many fibered knots for fixed genus $g \geq 2$.  

\begin{defn}\emph{(Dehn Surgery)} Let $K$ be a knot in the 3-sphere $S^{3}$. A Dehn surgery on $S^{3}$ is the surgery operation that removes tubular neighborhood of knot $K$ and gluing in the copy of solid torus $J = S^{1}\times D^{2}$ using the any diffeomorphism of the boundary torus. Let $\gamma = pm + ql$ be an essential loop on the boundary torus such that the meredian of the solid torus $J$ is glued to the curve $\gamma$, where $(m,l)$ is the standard meridian and longitude for knot $K$ and $p$, $q$ are coprime. The rational number $p/q$ or $\infty$, slope, is assigned to each surgery. \end{defn}

\subsection{Generalized fiber sum}
\label{sec:def}

\begin{defn} 
Let $X$\/ and $Y$\/ be closed oriented smooth\/ $4$-manifolds each containing a smoothly embedded surface $\Sigma$ of genus $g \geq 1$. Assume $\Sigma$ represents a homology of infinite order and has self-intersection zero in $X$\/ and $Y$, so that there exists a product tubular neighborhood, say $\nu\Sigma\cong \Sigma\times D^{2}$, in both $X$\/ and $Y$. Using an orientation-reversing and fiber-preserving diffeomorphism $ \psi : \Sigma\times S^1  \rightarrow \Sigma\times S^1$, we can glue $X \setminus\nu\Sigma$ and $Y\setminus \nu\Sigma$\/ along the boundary $\partial(\nu\Sigma)\cong \Sigma\times S^{1}$.  The resulting closed oriented smooth\/ $4$-manifold, denoted $X\#_{\psi}Y$, is called a \emph{generalized fiber sum}\/ of $X$\/ and $Y$\/ along $\Sigma$.  \end{defn}

\begin{defn}
Let $e(X)$ and $\sigma(X)$ denote the Euler characteristic and the signature
of a closed oriented smooth\/ $4$-manifold $X$, respectively.  We define 
\begin{equation*}
c_1^2(X):= 2e(X)+3\sigma(X), \quad \chi_h(X):=\frac{e(X)+\sigma(X)}{4}.
\end{equation*}
\end{defn}

\begin{lem}\label{lemma:c_1^2 and chi} 
Let $X$ and $Y$ be closed, oriented, smooth\/ $4$-manifolds containing an embedded surface\/ $\Sigma$ of self-intersection\/ $0$. Then 
\begin{eqnarray*}
c_{1}^{2}(X\#_{\psi}Y) &=&  c_{1}^{2}(X) + c_{1}^{2}(Y) + 8(g-1),\\
\chi_{h}(X\#_{\psi}Y) &=&  \chi_{h}(X) + \chi_{h}(Y) + (g-1),
\end{eqnarray*}
where $g$ is the genus of the surface $\Sigma$. 
\end{lem}

\begin{proof}
The above formulas simply follow from the well-known formulas 
\begin{equation*}
e(X\#_{\psi}Y)= e(X) + e(Y) - 2e(\Sigma),\quad 
\sigma(X\#_{\psi}Y) = \sigma(X) + \sigma(Y).  \qedhere   
\end{equation*}
\end{proof}
 
\medskip

If $X$\/ and $Y$\/ are symplectic $4$-manifolds and $\Sigma$ is a symplectic 
submanifold in both, then according to a theorem of Gompf \cite{Go}, 
$X\#_{\psi}Y$\/ admits a symplectic structure.  In such a case, we will call 
$X\#_{\psi}Y$\/ a \emph{symplectic sum}.

\section{Seiberg-Witten Invariants}

In this section, we review the basics of Seiberg-Witten invariants introduced by Seiberg and Witten. Let us recall that the Seiberg-Witten invariant  of a smooth closed oriented $4$-manifold $X$ with $b_2 ^+(X) > 1$ is an integer valued function which is defined on the set of $spin ^{\, c}$ structures over $X$ \cite{W}. For simplicity, we assume that $H_1(X,\Z)$ has no $2$-torsion. Then there is a natural identification of the $spin ^{\, c}$ structures of $X$ with the characteristic elements of $H^2(X,\Z)$ as follows: to each $spin ^c$ structure $\fs $ over $X$ corresponds a bundle of positive spinors $W^+_{\fs}$ over $X$. Let $c(\fs)\in H_2(X)$ denote the Poincar\'e dual of $c_1(W^+_{\fs})$. Each $c(\fs)$ is a characteristic element of $H_2(X,\Z)$ (i.e. its Poincar\'e dual $\hat{c}(\fs)=c_1(W^+_{\fs})$ reduces mod~2 to $w_2(X)$).

In this set up, we can view the Seiberg-Witten invariant as an integer valued function

\[ \ssw_X: \lbrace K\in H^2(X,\Z)|K\equiv w_2(TX)\pmod2\rbrace
\rightarrow \Z. \]  

\noindent If $\ssw_X(\b)\neq 0$, then we call $\b$ a {\it{basic class}} of $X$. It is a fundamental fact that the set of basic classes is finite. Furthermore, if $\b$ is a basic class, then so is $-\b$ with
\[\ssw_X(-\b)=(-1)^{(e+\sigma)(X)/4}\,\ssw_X(\b)\] where
$e(X)$ is the Euler characteristic and $\sigma(X)$ is the signature of $X$.

\smallskip

When $b_2 ^+(X) > 1$, then Seiberg-Witten invariant is a diffeomorphism invariant. It does not depend on the choice of generic metric on $X$ or a generic pertubation of Seiberg-Witten equations. 

\smallskip

If $b_2 ^+(X) = 1$, then the Seiberg-Witten invariant depends on the chosen metric and  perturbation of Seiberg-Witten equations. Let us recall that the perturbed Seiberg-Witten moduli space $\CM_{X}(\b, g, h)$ is defined as the solutions of the Seiberg-Witten equations 
           
\medskip
   
\begin{center} 

 $F_A^{+} = q(\psi) + ih, \ D_{A}\psi = 0$ 

\end{center}

\medskip

\noindent divided by the action of the gauge group. Where $A$ is connection on the line bundle $L$ with $c_{1}(L) = \b$, $g$ is riemannian metric on $X$, $\psi$ is the section of the positive spin bundle corresponding to the $spin^{c}$ structure determined by $\b$, $F_{A} ^{+}$ is the self-dual part of the curvature $F_{A}$, $D_{A}$ is the twisted Dirac operator, $q$ is a quadratic function, and $h$ is self-dual 2-form on $X$. If $b_{2}^{+}(X) \geq 1$ and $h$ is generic metric, then Seiberg-Witten moduli space $\CM_{X}(\b, g, h)$ is closed manifold with dimension $d = (\b^{2} -2e(X) - 3\sigma(X))/4$. The Seiberg-Witten invariant defined as follows: 

\smallskip

\begin{displaymath}
\left\{\begin{array}{ll}   
\ssw_{X} (\b) = \langle [ \CM_{X} (\b, g, h)], \mu^{d/2}\rangle & \textrm{if $d \geq 0$ and even}  \\
 \\
\ssw_{X} (\b) = 0 & \textrm{otherwise} \\
\end{array} \right.
\end{displaymath}

\noindent where $ \mu  \in  H^{2}(\CM_{X} (\b, g, h), \Z)$ is the Euler class of the base fibration.

When $b_{2}^{+}(X) = 1$, the Seiberg-Witten invariant $SW_{X}(\b, g, h)$ depends on $g$ and $h$. Because of this, there are two types of Seiberg-Witten invariants: $\ssw_{X}^{+}$ and $\ssw_{X}^{-}$.        

\smallskip

\begin{thm}\cite{KM}, \cite{OS} Let $X$ be closed, oriented, smooth 4-manifold with $b_{1}(X) = 0$ and $b_{2}^{+}(X) = 1$. Fix a homology orientation of $H_{+}^{2}(X, \R)$. For given riemannian metric $g$ let $\o_{+}^{g}$ be the unique $g$-harmonic self-dual 2-form that has norm $1$ and is compatible with the orientation of $H_{+}^{2}(X, \R)$. Then for each characteristic element $\b$ with $(\b^{2} -2e(X) - 3sign(X))/4 \geq 0$ the following holds: If $(2\pi \b + h_{1}) \ \cdot \ \o_{+}^{g_{1}}$ and $(2\pi \b + h_{2}) \ \cdot \ \o_{+}^{g_{2}}$ are not zero and have same sign, then $\ssw_{X}(\b, g_{1}, h_{1}) = \ssw_{X}(\b, g_{2}, h_{2})$.

\medskip

\end{thm}

\begin{defn} If $(2\pi \b + h) \ \cdot \ \o_{+}^{g} > 0$, then write $\ssw_{X}^{+}(\b)$ for $\ssw_{X}(\b, g, h)$. If $(2\pi \b + h) \ \cdot \  \o_{+}^{g} < 0$, then write $\ssw_{X}^{-}(\b)$ for $\ssw_{X}(\b, g, h)$  

\medskip

\end{defn}

\begin{thm}\cite{Sz} Let $X$ be closed, oriented, smooth 4-manifold with $b_{1}(X) = 0$ and $b_{2}^{+}(X) = 1$ and $b_{2}^{-} \leq 9$ . Then for each characteristic element $\b$, pairs of riemannian metrics $g_{1}$, $g_{2}$ and small perturbing 2-forms $h_{1}$, $h_{2}$ then $\ssw_{X}(\b, g_{1}, h_{1}) = \ssw_{X}(\b, g_{2}, h_{2})$.
  
\end{thm}

\begin{proof}. Let $\b$ be a characteristic element for which $d \geq 0$. Then $2e(X) + 3sign(X) = 4 + 5b_{2}^{+} - b_{2}^{-} \geq 0$, which as implies $\b^{2} \geq 0$. As a corollary, it follows that $\b \cdot \ \o_{+}^{g_{1}}$ and $\b \cdot \  \o_{+}^{g_{2}}$ are both non-zero and have same signs. Now using the Theorem 3.1, the result follows.

\end{proof} 

\smallskip

\begin{thm}\cite{LL} {\bf (Wall crossing formula)} Assume that $X$ is closed, oriented, smooth 4-manifold with $b_{1}(X) = 0$ and $b_{2}^{+}(X) = 1$. Then for each characteristic line bundle $L$ on $X$ such that the formal dimension of the Seiberg-Witten moduli space is non-negative even integer $2m$, then $\ssw_{X}^{+}(L) - \ssw_{X}^{-}(L) = -(-1)^{m}$.

\end{thm}

Note that when $b_{2}^{-} \leq 9$, it follows from the above result that there is well defined Seiberg-Witten invariant which will be denoted as $\ssw_{X}^{o}(X)$. 

\medskip

\begin{thm}\cite{T} Suppose that $X$ is closed symplectic 4-manifold with $b_{2}^{+}(X) > 1$ $(b_{2}^{+}(X) = 1)$. If $K_{X}$ is a canonical class of $X$, then $\ssw_{X}(\pm K_{X}) = \pm 1$ $(SW_{X}^{-}(-K_{X}) = \pm 1)$. \end{thm}

\medskip

\begin{defn} The 4-manifold $X$ is of simple type if each basic class $\b$ satisfies the equation ${\b}^2 = {c_{1}}^{2}(X) = 3\sigma(X) + 2e(X)$. If $X$ is symplectic manifold of $b_{2}^{+}(X) > 1$ then it has simple type.  

\end{defn}

\medskip

\begin{thm}\cite{KM}, \cite{OS} {\bf (Generalized adjunction formula for $b_{2}^{+} > 1$)} Assume that $\Sig \subset X$ is an embedded, oriented, connected surface of genus $g(\Sig)$ with self-intersection $|\Sig|^{2} \geq 0$ and represents nontrivial homology class. Then for every Seiberg-Witten basic class $\b$, $2g(\Sig) - 2 \geq |\Sig|^{2} + |\b(|\Sig|)|$. If $X$ is of simple type and $g(\Sig) > 0$, then the same inequality holds for $\Sig$ with arbitrary square $|\Sig|^{2}$.

\end{thm}   

\medskip

\begin{thm}\cite{LL} {\bf (Generalized adjunction formula for $b_{2}^{+} = 1$)} Let $\Sig \subset X$ is an embedded, oriented, connected surface of genus $g(\Sig)$ with self-intersection $|\Sig|^{2} \geq 0$ and represents nontrivial homology class. Then any characteristic class $\b$ with ${\ssw_{X}}^{o} (\b) \ne 0$ satisfy the following inequality $2g(\Sig) - 2 \geq \Sig^{2} + |\b(|\Sig|)|$.

\end{thm}

\section {Symplectic manifolds cohomology equivalent to $S^{2}\times S^{2}$}

To construct our manifolds, we start with symplectic $4$-manifolds described below. By applying Gompf's symplectic fiber sum operation, we will obtain the manifolds $X_{K}$.

Let $K$ be a genus one fibered knot (i.e. trefoil or figure eight knot) in $S^{3}$ and $m$ a meridional circle to $K$. Perform 0-framed Dehn surgery on $K$ and denote the resulting $3$-manifold by $M_{K}$. The manifold $M_{K}$ has same integral homology as $S^{2}\times S^{1}$, where class of $m$ generates $H_{1}(M_{K})$. Since the knot $K$ is genus one fibered knot, it follows that the manifold $M_{K}\times S^{1}$ is a torus bundle over the torus which is homology equivalent to $T^{2}\times S^{2}$. Since $K$ is a fibered knot, $M_{K}\times S^{1}$ admits a symplectic structure. Note that $m \times S^{1}$ is the section of this fibration. The first homology of $M_{K}\times S^{1}$ is generated by the standard first homology classes of the torus section and the classes $\gamma_{1}$ and $\gamma_{2}$ of the fiber $F$ of the given fibration is trivial in the homology. Smoothly embedded torus section $T_{m} = m\times S^{1}$ has a self-intersection zero and its neighborhood in $M_{K}\times S^{1}$ has a canonical identification with $T_{m}\times D^{2}$.

The intermediate building block in our construction will be the twisted fiber sum of the two copies of manifold $M_{K}\times S^{1}$,  where we identify the fiber of one fibration to the base of other. Let $Y_{K}$ denote the mentioned twisted fiber sum $Y_{K} = M_{K}\times S^{1} \#_{F = T_{m}}\,\ M_{K}\times S^{1}$. It follows from Gompf's theorem \cite {GS} that $Y_{K}$ is symplectic.  

Let $T_{1}$ be the section of the first copy of $M_{K}\times S^{1}$ in $Y_{K}$ and $T_{2}$ be the fiber of the second copy. Then $T_{1} \# T_{2}$ symplectically embeds into $Y_{K}$ \cite {FS2}. Now suppose that $Y_{K}$ is the symplectic $4$-manifold given, and $\Sigma_{2} = T_{1} \# T_{2}$ is the genus two symplectic submanifold of self-intersection zero sitting inside of $Y_{K}$. Let $(m, x, \gamma_{1},\gamma_{2})$ be the generators of $H_{1}(\Sigma_{2})$ under the inclusion-induced homomorphism. We choose the diffeomorphism  $\phi : T_{1} \# T_{2} \longrightarrow T_{1} \# T_{2}$ of the $\Sigma_{2}$ that changes the generators of the first homology according to the following rule: $\phi (m') = \gamma_{1}$, $\phi ({\gamma_1}') = m$, $\phi (x') = \gamma_{2}$ and $\phi ({\gamma_2}') = x$. Next, we take the fiber sum along the genus two surface $\Sig_{2}$ and denote the new symplectic manifold as $X_{K}$ i.e. $X_{K} = Y_{K}\#_{\phi}\,\ Y_{K}$. We will show that the new manifold $X_{K}$ has a trivial first Betti number and has same integral cohomology as $S^{2}\times S^{2}$. Furthermore, $e(X_{K}) = 4$, $\sigma (X_{K}) = 0$, ${c_1}^{2}(X_{K}) = 8$, and $\chi_{h}(X_{K}) = 1$. We will compute $H_{1}(X_{K})$ first by applying Mayer-Vietoris sequence and next by directly computing the fundamental group of $X_{K}$.

\begin{lem} $H_{1} (X_{K}, \Z) = 0$ and $H_{2} (X_{K}, \Z) = \Z \oplus \Z$.

\end{lem}

\begin{proof} We use Mayer-Vietoris sequence to compute the homology of $X_{K} = Y_{K}\#_{\phi}\,\ Y_{K}$. Let $Y_{1} = Y_{2} = Y_{K} \setminus \nu {\Sigma_2}$. Note that $Y_{1} \cap Y_{2}$ is homologous to $\Sigma_{2} \times S^{1}$. By applying the reduced Mayer-Vieotoris to the triple $(X_{K}, Y_{1}, Y_{2})$, we have the following long exact sequence \\

$\cdots \longrightarrow H_{2}(S^{1} \times \Sig_{2}) \longrightarrow H_{2}(Y_{1}) \oplus H_{2}(Y_{2})  \longrightarrow H_{2}(X_{K})  \longrightarrow H_{1}(\Sig_{2}\times S^{1})  \longrightarrow H_{1}(Y_{1}) \oplus H_{1}(Y_{2})  \longrightarrow H_{1}(X_{K}) \longrightarrow 0$  \\

The simple computation by Kunneth formula yields $H_{1} (\Sig_{2}\times S^{1}, \Z) = \Z \oplus \Z \oplus \Z \oplus \Z \oplus \Z = \ <\lambda> \oplus <m> \oplus <x> \oplus<\gamma_{1}> \oplus< \gamma_{2}>$. Also, we have $H_{1} (Y_{1}, \Z) = \Z \oplus \Z = \ <m> \oplus <x>$ and  $H_{1} (Y_{2}, \Z) = \Z \oplus \Z = \ <m'> \oplus <x'>$ \\

Let $\phi_{*}$ and $\delta$ denote the last two arrows in the long exact sequence above. Because the way map $\phi$ is defined, the essential homology generators will map to trivial ones, thus $\phi_{*}(m) = \phi_{*}(x) = \phi_{*}(m') = \phi_{*}(x')$ = $0$. Since $Im(\phi_{*}) = Ker(\delta)$, we conclude that $H_{1}(X_{K}) = Ker(\delta) = Im(\phi_{*}) = 0$  

Next, by using the facts that $b_{1} = b_{3} = 0$, $b_{0} = b_{4} = 1$ and symplectic sum formula for Euler characteristics, we compute $b_{2} = e(Y_{K}) + e(Y_{K}) + 2 = 0 + 0 + 2 = 2$

We conclude that $H_{2}(X_{K}, \Z) = \Z \oplus \Z$. Note that a basis for the second homology are classes of $\Sigma_{2} = S$ and the new genus two surface $T$ resulting from the last fiber sum opeartion, where $S^{2} = T^{2} = 0$ and $S \cdot T = 1$. Thus, the manifolds obtained by the above construction have intersection form $H$, so they are spin. \end{proof}

\medskip

\begin{lem} $e(X_{K}) = 4$, $\sigma (X_{K}) = 0$, ${c_1}^{2}(X_{K}) = 8$, and $\chi_{h}(X_{K}) = 1$.

\end{lem}

\begin{proof} Using the lemma 2.10, we have $e(X_{K}) = 2e(Y_{K}) + 4$, $\sigma (X_{K}) = 2\sigma(Y_{K})$, ${c_1}^{2}(X_{K}) = 2{c_1}^{2}(Y_{K}) + 8$, and $\chi_{h}(X_{K}) = 2\chi_h(Y_{K}) + 1$ . Since $e(Y_{K}) = 0$, $\sigma (Y_{K}) = 0$, ${c_1}^{2}(Y_{K}) = 0$ and $\chi_{h}(Y_{K}) = 0$ result follows. \end{proof}

\subsection{Fundamental Group Computation for Trefoil}

\subsubsection{Step 1: Fundamental Group of $M_{K}\times S^{1}$} In this section we will assume that $K$ is trefoil. The case when $K$ is figure eight can be treated similarly.  

Let $a$ and $b$ denote the Wirtinger generators of the trefoil knot. Then the group of $K$ has the following presentation 

\begin{equation*}
\pi_1(S^3\setminus \nu K) \ = \langle a, b \mid  aba=bab \rangle \ = \langle u, v \mid  u^{2} =  v^{3} \rangle \end{equation*}

\noindent where $u = bab$ and $v = ab$. By Lemma 2.2, the homotopy classes of the meridian and the longitude of the trefoil are given as follows: $ m = uv^{-1} = b$ and $l = u^{2}(uv^{-1})^{-6} = ab^2ab^{-4}$. Notice that $\gamma_{1} = a^{-1}b$ and $\gamma_{2} = b^{-1}aba^{-1}$ generate the image of the fundamental group of the Seifert surface of $K$ under the inclusion-induced homomorphism. 
 Let $M_{K}$ denote the result of 0-surgery on $K$. 

\begin{lem} \begin{eqnarray*}
\pi_{1}(M_{K}\times S^1) &=& \pi_1(M_K) \oplus \Z \\ &=&  \langle a,b,x \mid aba = bab,\, ab^{2}ab^{-4} = 1,\, [a,x]=[b,x]=1 \rangle.
\end{eqnarray*}

\end{lem}

\begin{proof}

Notice that the fundamental group of $M_{K}$ is obtained from the knot group of the trefoil by adjoining the relation $l = u^{2}(uv^{-1})^{-6} = ab^2ab^{-4} = 1$. Thus, we have the presentation given above.

\end{proof}

\subsubsection{Step 2: Fundamental Group of $Y_{K}$} Next, we take the two copies of the manifold $M_{K}\times S^{1}$. In the first copy, take a tubular neighborhood of the torus section $b\times x$, remove it from $M_{K}\times S^{1}$ and denote the resulting manifold as $C_{S}$. In the second copy, we remove a tubular neighborhood of the fiber $F$ and denote the complement by $C_{F}$. 

\medskip

\begin{lem}
Let $C_S$ be the complement of a neighborhood of a section in $M_K\times S^1$.
Then we have 
\begin{equation*}
\pi_1(C_S) = \langle 
a, b, x \mid  aba=bab,\, [a,x]=[b,x]= 1
\rangle.  
\end{equation*}  
\end{lem}

\begin{proof}
Note that $C_S=(M_K\setminus \nu (b))\times S^1=(S^3\setminus \nu K)
\times S^1$.  
\end{proof}

\begin{lem}
Let $C_F$ be the complement of a neighborhood of a fiber in $M_K\times S^1$.  Then we have 
\begin{eqnarray*}
\pi_1(C_F)&=&\langle \gamma_{1}',\gamma_{2}', d, y \mid
[\gamma_{1}',\gamma_{2}']=[y,\gamma_{1}']=[y,\gamma_{2}']=1,\\ 
&&d\gamma_{1}'d^{-1}=\gamma_{1}'\gamma_{2}',\,
d\gamma_{2}'d^{-1}=(\gamma_{1}')^{-1} \rangle.
\end{eqnarray*}
\end{lem}

\begin{proof}
To compute the fundamental group of $C_{F}$, we will use the following observation: $C_F$ is homotopy equivalent to a torus bundle over a wedge of two circles. Thus the generators $d$ and $y$ do not commute in the fundamental group of $C_{F}$. Also, the monodromy along the circle $y$ is trivial whereas the monodromy along the circle $d$ is the same as the monodromy of $M_K$. This implies that the fundamental group of $C_{F}$ has the presentation given above.  \end{proof}

\begin{lem}
Let $Y_K$ be the symplectic sum of two copies of $M_K\times S^1$,
identifying a section in one copy with a fiber in the other copy.  If the gluing map
$\psi$ satisfies $\psi_{\ast}(x)=\gamma'_1$ and $\psi_{\ast}(b)=\gamma'_2$, then 
\begin{eqnarray*}
\pi_1(Y_K)&=&\langle 
a,b,x, \gamma_1',\gamma_2',d,y
\mid
aba=bab,\, [x,a]=[x,b]=1,\\
&& [\gamma_1',\gamma_2']=[y,\gamma_1']=[y,\gamma_2']=1,\, 
d\gamma_1'd^{-1}=\gamma_1'\gamma_2',\,
d\gamma_2'd^{-1}=(\gamma_1')^{-1},\\
&&x=\gamma_1',\, b=\gamma_2',\, ab^2ab^{-4} =[d,y]
\rangle \\
&=&\langle a,b,x,d,y
\mid aba=bab,\, [x,a]=[x,b]=1,\\
&&[y,x]=[y,b]=1,\, dxd^{-1}=xb,\, dbd^{-1}=x^{-1},\,
ab^2ab^{-4}=[d,y] \rangle.
\end{eqnarray*}
\end{lem}

\begin{proof}
By Van Kampen's Theorem, $\pi_1(Y_K)=\pi_1(C_S)\ast\pi_1(C_F)/
\pi_1(T^3)$.  One circle factor of $T^3$ is identified with the longitude
of $K$ on one side and the meridian of the torus fiber in $M_K\times S^1$
on the other side.  This gives the last relation.  
\end{proof}

Inside $Y_K$, we can find a genus 2 symplectic submanifold 
$\Sigma_2$ which is the internal sum of a punctured fiber in $C_S$
and a punctured section in $C_F$.  The inclusion-induced homomorphism
maps the standard generators of $\pi_1(\Sigma_2)$ to $a^{-1}b$, 
$b^{-1}aba^{-1}$, $d$\/ and $y$.  

\begin{lem}\label{lemma:Y_K complement}
There are nonnegative integers $m$ and $n$ such that
\begin{eqnarray}
\pi_1(Y_K\setminus \nu \Sigma_2 )&=&\langle
a,b,x,d,y;\, g_1,\dots,g_m \mid aba=bab,\, \label{eq:Y_K complement} \\
&&[y,x]=[y,b]=1,\, dxd^{-1}=xb,\, dbd^{-1}=x^{-1},\, \nonumber\\
&&ab^2ab^{-4}=[d,y], \,
r_1=\cdots=r_n=1,\, r_{n+1}=1 \rangle, \nonumber
\end{eqnarray}
where the generators $g_1,\dots,g_m$ and relators
$r_1,\dots,r_n$ all lie in the normal subgroup $N$\/ 
generated by the element\/ $[x,b]$, and the relator $r_{n+1}$
is a word in $x,a$ and elements of $N$.
Moreover, if we add an extra relation\/ 
$[x,b]=1$ to\/ $(\ref{eq:Y_K complement})$, then the relation\/ 
$r_{n+1}=1$ simplifies to\/ $[x,a]=1$.
\end{lem}

\begin{proof}
This follows from Van Kampen's Theorem.  Note that 
$[x,b]$ is a meridian of $\Sigma_2$ in $Y_K$.  Hence setting $[x,b]=1$
should turn $\pi_1(Y_K\setminus \nu \Sigma_2 )$ 
into $\pi_1(Y_K)$.   Also note that $[x,a]$ is the boundary of 
a punctured section 
in $C_S\setminus \nu ({\rm fiber})$, and is no longer trivial in
$\pi_1(Y_K\setminus\nu\Sigma_2)$.  By setting
$[x,b]=1$, the relation $r_{n+1}=1$ is to turn into $[x,a]=1$.

It remains to check that the relations in $\pi_1(Y_K)$
other than $[x,a]=[x,b]=1$ remain the same in $\pi_1(Y_K\setminus\nu\Sigma_2)$.  
By choosing a suitable point $\theta\in S^1$ away from the image of
the fiber that forms part of $\Sigma_2$, we obtain an embedding 
of the knot complement 
$(S^3\setminus\nu K)\times\{\theta\} \hookrightarrow
C_S\setminus \nu ({\rm fiber})$.  This means that $aba=bab$\/ holds in 
$\pi_1(Y_K\setminus\nu\Sigma_2)$.   Since $[\Sigma_2]^2=0$, there 
exists a parallel copy of $\Sigma_2$ outside $\nu\Sigma_2$, wherein
the identity $ab^2ab^{-4}=[d,y]$ still holds.  The other remaining relations in 
$\pi_1(Y_K)$ are coming from the monodromy of the torus bundle over a
torus.  Since these relations will now describe the monodromy of the 
punctured torus bundle over a punctured torus,
they hold true in $\pi_1(Y_K\setminus\nu\Sigma_2)$.    
\end{proof}

\subsubsection {Step 3: Fundamental Group of $X_{K}$} 

Now take two copies of $Y_K$.  Suppose that the fundamental group of the second 
copy has $e,f,z,s,t$\/ as generators, and 
the inclusion-induced homomorphism in
the second copy maps the generators of $\pi_1(\Sigma_2)$ to $e^{-1}f$, 
$f^{-1}efe^{-1}$, $s$\/ and $t$.  Let $X_K$ denote the symplectic sum of two copies
of $Y_K$ along $\Sigma_2$, where the gluing map $\psi$ maps the generators as 
follows:
\begin{equation*}
\psi_{\ast}(a^{-1}b)=s,\; 
\psi_{\ast}(b^{-1}aba^{-1})= t,\; 
\psi_{\ast}(d)= e^{-1}f,\; 
\psi_{\ast}(y)=  f^{-1}efe^{-1}.
\end{equation*}

\begin{lem}\label{lemma:pi_1 X_K}
There are nonnegative integers $m$ and $n$ such that 
\begin{eqnarray*}
\pi_1(X_K)&=&\langle a,b,x,d,y;\;
e,f,z, s, t;\; g_1,\dots, g_m;\; h_1,\dots, h_m \mid \\
&&aba=bab,\, [y,x]=[y,b]=1,\\
&& dxd^{-1}=xb,\, dbd^{-1}=x^{-1},\,
ab^2ab^{-4}=[d,y],\\
&&r_1=\cdots=r_n=r_{n+1}=1,\, r_1'=\cdots=r_n'=r'_{n+1}=1 ,\, \\
&&efe=fef,\, [t,z]=[t,f]=1,\,\\
&& szs^{-1}=zf,\, sfs^{-1}=z^{-1},\,
ef^2ef^{-4}=[s,t],\\
&& d=e^{-1}f,\, y=f^{-1}efe^{-1},\, 
a^{-1}b=s,\, b^{-1}aba^{-1}=t,
[x,b]=[z,f] \rangle ,
\end{eqnarray*}
where $g_i,h_i$ $(i=1,\dots,m)$ and $r_j,r_j'$ $(j=1,\dots,n)$ all 
lie in the normal subgroup $N$ generated by
$[x,b]=[z,f]$.  Moreover, $r_{n+1}$ is a word in $x,a$ and elements of $N$, 
and $r'_{n+1}$ is a word in $z,e$ and elements of $N$.   
\end{lem}

\begin{proof}
This is just a straightforward application of Van Kampen's Theorem and
Lemma \ref{lemma:Y_K complement}. Also, notice that the abelianization of $\pi_1(X_K)$ is trivial.
\end{proof}

\noindent The different symplectic cohomology $S^{2}\times S^{2}$'s can be obtained if we use other genus $1$ fibered knot, the figure eight knot, or combination of both in our construction. These manifolds can be distinguished by their fundamental groups. In addition, using a family of non fibered genus one twist knots, we also obtain an infinite family of cohomology $S^{2}\times S^{2}$. This family of cohomology $S^{2}\times S^{2}$'s will not be symplectic anymore and can be distinguished by Seiberg-Witten invariants.

\subsection{Seiberg-Witten invariants for manifold $X_{K}$}

Let $C$ be a basic class of the manifold $X_{K}$. We can write $C$ as a linear combination of $S$ and $T$, i.e. $C = aS + bT$. $X_{K}$ is a symplectic manifold and has simple type. So for any basic class $C$, $C^{2} = 3 \sigma(X_{K}) + 2e(X_{K}) = 8$. It follows that $2ab = 8$. Next we apply the adjunction inequality to $S$ and $T$ to get $2g(S) - 2 \geq  [S]^2 + |C(S)|$ and $2g(T) - 2 \geq  [T]^2 + |C(T)|$. These gives us two more restriction on $a$ and $b$ : $2 \geq |b|$ and $2 \geq |a|$. Thus, it follows that $C = \pm (2S + 2T)$, which is $\pm$ the canonical class $K_{X_{K}} = 2S + 2T$ of $X_{K}$. Notice that, since $b_{2}^{-}(X_{K}) \leq 9$ and $(-K_{X_{K}})\cdot \o <0$, we have well defined $\ssw_{X_{K}}^{o}$. Now it follows from the theorems of Section 2 that $\ssw_{X_{K}}^{0}(-K_{X_{K}} ) = \ssw_{X_{K}}^{-}( -K_{X_{K}} ) = \pm 1$.

\section {Symplectic manifolds cohomology equivalent to $\#_{(2g-1)}{(S^{2}\times S^{2})}$}

We will modify our construction in the previous section to get an infinite family of symplectic 4-manifolds cohomology equivalent to $\#_{(2g-1)}(S^{2}\times S^{2})$  for any $g \geq 2$. 
  
Let $K'$ denote the genus $g$ fibered knot in $S^{3}$ and $m$ a meridional circle to $K'$. We first perform $0$-framed surgery on $K'$ and denote the resulting $3$-manifold by $Z_{K'}$. The $4$-manifold $Z_{K'}\times S^{1}$ is $\Sigma_{g}$ bundle over the torus and has the same integral homology as $T^{2}\times S^{2}$. Since $K'$ is a fibered knot, $Z_{K'}\times S^{1}$ is symplectic manifold. Again, there is a torus section $m \times S^{1} = T_{m}$ of this fibration. The first homology of $Z_{K'}\times S^{1}$ is generated by the standard first  homology classes of this torus section. The standard homology generators of the fiber $F$, which we denote as $\alpha_{2}, \beta_{2}, \cdots, \alpha_{g+1}$,  and $\beta_{g+1}$ of the given bundle is trivial in the homology. The section $T_{m}$ has zero self-intersection and the its neighborhood in $Z_{K'}\times S^{1}$ has a canonical identification with $T_{m}\times D^{2}$.

We form the twisted fiber sum of the manifold $M_{K}\times S^{1}$ with $Z_{K'}\times S^{1}$ where we identify fiber of the first fibration to the section of other. Let $W_{K'}$ denote the new manifold $W_{K'} = M_{K}\times S^{1} \#_{F = T_{m'}}\,\ Z_{K'}\times S^{1}$. Again, it follows from Gompf's theorem \cite {Go} that $W_{K'}$ is symplectic. 

Let $T_{1}$ be the section of the $M_{K}\times S^{1}$ which we discussed earlier and $\Sigma_{g} = \Sig^{'}$ be the fiber of the $Z_{K'}\times S^{1}$. Then $T_{1} \# \Sig_{g}$ embeds in $W_{K'}$ and has self-intersection zero. Now suppose that $W_{K'}$ is the symplectic 4-manifold given, and $\Sigma_{g+1} = T_{1} \# \Sig_{g}$ is the genus $g+1$ symplectic submanifold of self-intersection $0$ sitting inside of $W_{K'}$. Let $\alpha_1$, $\beta_1$, $\alpha_2$, $\beta_2$, $\dots$ , $\alpha_{g+1}$, and $\beta_{g+1}$ be the generators of the first homology of the surface $T_{1} \# \Sig_{g} = \Sig_{g+1}$. Let $\psi : T^{2} \# \Sig_{g} \longrightarrow \Sig_{g} \# T^{2}$ be the diffeomorphism of the genus $g+1$ surface that changes the generators of the first homology according to the following rule: $\psi (\alpha_1) = \alpha_{g+1}$, $\psi (\alpha_{g+1}) = \alpha_{1}$, $\psi (\beta_{1}) = \beta_{g+1}$, and $\psi (\beta_{g+1}) = \beta_{1}$. Take the fiber sum along this genus $g+1$ surface $\Sigma_{g+1}$ and denote the resulting symplectic manifold as $V_{K'}$. The new manifold $V_{K'} = W_{K'}\#_{\psi}\,\ W_{K'}$ has trivial first homology and has same homology of $\#_{(2g-1)}S^{2}\times S^{2}$. 

\begin{lem} $H_{1} (V_{K'}, \Z) = 0$, $H_{2} (V_{K'}, \Z) = \Z \oplus \Z \oplus \Z \oplus \Z \cdots \oplus \Z$ where there are $2(2g-1)$ copies of $\Z$

\end{lem}

\begin{proof}

The proof is similar to genus one case and can be obtained by applying Mayer-Vietoris sequence.

\end{proof}

\noindent Note that $H_{2}(V_{K'}, \Z)$ has rank $4g - 2$. The base for second homology are classes of the $\Sig_{g+1} = S$, the genus two surface $\Sigma_{2} = \Sig$  resulting from the last fiber sum operation, $(2g-2)$ rim toris $R_{i}$ and $(2g-2)$ assosiated vanishing classes $V_{i}$ in the $\Sig^{'}$ \cite{FS2}. In the manifold $V_{K'}$ these classes contribute $2g-2$ new hyperbolic pairs. Thus, the manifolds obtained by the above construction have intersection form $\oplus_{2g-1}H$. Notice that $b_{2}^{+}(V_{K'}) = b_{2}^{-}(V_{K'}) = 2g-1$

\begin {lem} $e (V_{K'}) = 4g$, $\sigma (V_{K'}) = 0$, ${c_1}^{2}(V_{K'}) = 8g$, and $\chi_{h}(V_{K'}) = g$

\end {lem}

\begin{proof} Using the lemma 2.9, we have $e(V_{K'}) = 2e(W_{K'}) + 4g$, $\sigma (V_{K'}) = 2\sigma(W_{K'})$, ${c_1}^{2}(V_{K'}) = 2{c_1}^{2}(W_{K'}) + 8g$ and $\chi_{h} (V_{K'}) = 2\chi_{h}(W_{K'}) + g$. Since $e(W_{K'}) = 0$ and $\sigma (W_{K'}) = 0$, our result follows.
\end{proof}

\subsection{Seiberg-Witten invariants for manifold $V_{K'}$}

Let $C$ be a basic class of the manifold $V_{K'}$. Let us write $C$ as a linear combination of $S$, $\Sigma$, rim toris $R_{i}$, $i=1,...,2g-2$ and the assosiated vanishing classes $V_{i}$, $i =1,\cdots ,2g-2$, $C = aS + b\Sig + \sum_{i=1}^{2g-2} u_{i}R_{i} + v_{i}V_{i}$. By adjunction inequality, the intersection number of any basic class with rim torus $0$ i.e. $C \cdot  R_{i} = 0$, $i=1,\cdots ,2g-2$. It follows that $Q^{T}v  = 0$ where $Q$ is the intersection matrix and $v = (v_{1},\cdots ,v_{2g-2})$. Using the fact that $Q$ is invertible, we get $v_{1} = ... = v_{2g-2} = 0$. Next, by applying the adjunction inequality again, we  have $V_{i} \cdot C = 0$ for $i = 1,...,2g-2$. Later give rise to system $Qu = 0$, which implies $u_{1} = \cdots = u_{2g-2} = 0$. This shows that any basic class has form $C = aS + b\Sigma$. Since $V_{K'}$ is a symplectic manifold, it has simple type. So for any basic class $C$, $C^{2} = 3 \sigma(V_{K'}) + 2e(V_{K'}) = 8g$. It follows that $2ab = 8g$. Next we apply the adjunction inequality to $S$ and $\Sigma$ to get $2g(S) - 2 \geq  [S]^2 + |C(S)|$ and $2g(\Sig) - 2 \geq  [\Sig]^2 + |C(\Sig)|$. These gives us two more restriction on $a$ and $b$ : $2g \geq |b|$ and $2 \geq |a|$. It follows that $C = \pm (2S + 2g\Sig)$, which are $\pm$ the canonical class of the given manifold. By applying the Taubes theorem from Section 2, we conclude that the value of Seiberg-Witten invariants on these classes is $\pm 1$.

\end{document}